\documentclass[10pt,a4paper]{paper}
\usepackage[english]{babel}

\usepackage{amsmath,amsfonts,amssymb,amsthm}
\usepackage{graphicx}
\usepackage{url}

\newtheorem{thm}{Theorem}
\newtheorem{prop}{Proposition}
\newtheorem{remark}{Remark}

\begin{document}

\title{Dengue disease, basic reproduction 
number and control\thanks{This is a preprint 
of a paper whose final and definitive form has appeared 
in \emph{International Journal of Computer Mathematics} (2011), 
DOI: 10.1080/00207160.2011.554540}}

\author{Helena Sofia Rodrigues$^1$\\
\url{sofiarodrigues@esce.ipvc.pt}
\and
M. Teresa T. Monteiro$^2$\\
\url{tm@dps.uminho.pt}
\and
Delfim F. M. Torres$^3$\\
\url{delfim@ua.pt}
\and
Alan Zinober$^4$\\
\url{a.zinober@sheffield.ac.uk}}

\institution{$^1$School of Business Studies,
Viana do Castelo Polytechnic Institute, Portugal\\[0.3cm]
$^2$Department of Production and Systems,
University of Minho, Portugal\\[0.3cm]
$^3$Department of Mathematics,
University of Aveiro, Portugal\\[0.3cm]
$^4$Department of Applied Mathematics,
University of Sheffield, UK}

\maketitle


\begin{abstract}
Dengue is one of the major international public health concerns.
Although progress is underway, developing a vaccine against the disease
is challenging. Thus, the main approach to fight the disease is vector control.
A model for the transmission of Dengue disease is presented.
It consists of eight mutually exclusive compartments representing the
human and vector dynamics. It also includes a control parameter (insecticide)
in order to fight the mosquito. The model presents three possible equilibria:
two disease-free equilibria (DFE) and another endemic equilibrium.
It has been proved that a DFE is locally asymptotically stable,
whenever a certain epidemiological threshold, known as the \emph{basic reproduction number},
is less than one. We show that if we apply a minimum level of insecticide,
it is possible to maintain the basic reproduction number below unity.
A case study, using data of the outbreak that occurred in 2009 in Cape Verde, is presented.

\medskip


\medskip

\end{abstract}


\section{Introduction}

The first recognized Dengue epidemic occurred almost
simultaneously in Asia, Africa and North America in the 1780s,
shortly after the identification and naming of the disease in
1779. It has spread especially in the tropical and subtropical
regions around the world, and nowadays is a disease widely found
in urban and semi-urban areas.

According to the World Health Organization, 50 to 100 million
Dengue infections occur yearly, including 500,000 Dengue
Haemorrhagic Fever cases and 22000 deaths, mostly among children
\cite{WHO}. Growing awareness of global climate change has
stimulated several assessments of its likely effects on
vector-borne disease as well as on health outcomes. Some of these
studies have indicated that countries with a mild climate, such as
in the Mediterranean, are at risk due to future climate conditions
that may be favourable to this kind of disease
\cite{Hopp2001}. This risk may be aggravated further due
to the volume of international tourism and trade that this region
experiences \cite{Semenza2009,SOL2009}. Travellers play an
essential role in the global epidemiology. They act as
viraemic travellers, carrying the disease into areas with mosquitos
that can transmit the infection.
This is particularly true with respect to the
reality in the archipelago of Cape Verde, where it is believed
that the insects responsible for the outbreak that occurred in 2009
came from Brazil, transported by means of the frequent air transport \cite{Africa21}.

There are two forms of Dengue: Dengue Fever and Dengue Haemorrhagic Fever.
The first one is characterized by a sudden fever without respiratory symptoms,
accompanied by intense headaches and lasts between three and seven days.
The second has the previous symptoms but additionally also nausea,
vomiting and fainting due to low blood pressure, and can lead to death
in two or three days \cite{Derouich2003}.

There are four distinct, but closely related, viruses that cause
Dengue. Recovery from infection by one virus provides lifelong
immunity against that virus but confers only partial and transient
protection against subsequent infection by the other three
viruses. There is good evidence that sequential infection
increases the risk of developing Dengue Haemorrhagic Fever (DHF).
Activities, such as triage and management, are critical in
determining the clinical outcome of Dengue. A rapid and efficient
front-line response not only reduces the number of unnecessary
hospital admissions but also saves lives.

There is no vaccine to protect against Dengue. With four closely
related viruses that can cause the disease, a vaccine would need
to immunize against all four types to be effective. There is
limited understanding of how the disease typically behaves and how
the virus interacts with the immune system. Another
challenge is that some studies show that some secondary Dengue
infection can leave to DHF, and theoretically a vaccine could be a
potential cause of severe disease if a solid immunity is not
established against the four serotypes. Research to develop
a vaccine is ongoing and the incentives to study the
mechanism of protective immunity are gaining more support, now that
the number of outbreaks around the world is increasing
\cite{Who2009}.

The spread of Dengue is attributed to expanding geographic
distribution of the four Dengue viruses and their mosquito
vectors, the most important of which is the predominantly urban
species \emph{Aedes aegypti}; some notes about it are presented in
Section~\ref{sec:2}. A mathematical model of the interaction
between human and mosquito populations is presented in
Section~\ref{sec:3}. Section~\ref{sec:4} is concerned with the
basic reproduction number, the equilibria of the epidemiological
model and their stability. In Section~\ref{sec:5}, previous
results are applied to a case study. Finally, concluding notes are
given in Section~\ref{sec:6}.


\section{Biological notes on \emph{Aedes aegypti}}
\label{sec:2}

The mosquito \emph{Aedes aegypti} is a tropical and
subtropical specie widely distributed around the world, mostly
between latitudes $35^{o}$N and 35$^o$S, which corresponds,
approximately, to a winter isotherm of 10$^o$C \cite{Braga2007}.
The life cycle of a mosquito has four distinct stages: egg, larva,
pupa and adult. In the case of \emph{Aedes aegypti}, the first
three stages take place in or near water whilst air is the medium
for the adult stage \cite{Otero2008}. The eggs of \emph{Aedes
aegypti} can resist desiccation and low temperatures for up to one
year. Although the hatching of mature eggs may occur spontaneously
at any time, this is greatly stimulated by flooding. The larva
moults four times in a period of a few days which culminates in
the pupal stage. The pupal stage lasts from one day to a few
weeks, depending on the temperature. At the end of this stage, the
adult emerges from the pupal skin \cite{Christophers1960}.
Studies suggest that most female mosquitoes may spend their
lifetime in or around the houses where they emerge as adults. This
means that people, rather than mosquitoes, rapidly move the virus
within and between communities.

The adult stage of the mosquito is considered to last an average
of eleven days in the urban environment. Dengue is spread only by
adult females that require a blood meal for the development of
eggs; male mosquitoes feed on nectar and other sources of sugar.
In this process, the female acquires the virus while feeding on the
blood of an infected person. After virus incubation for eight to
ten days, an infected mosquito is capable, during probing and
blood feeding, of transmitting the virus for the rest of its life.
\emph{Aedes aegypti} is one of the most efficient vectors
for arboviruses because it is highly anthropophilic, frequently
bites several times before complete oogenesis and thrives in
close proximity to humans \cite{Who2009}.

It is very difficult to control or eliminate the \emph{Aedes
aegypti} mosquito because it adapts to the environment and becomes
resistant to natural phenomena, \textrm{e.g.}, droughts, or human
interventions and control measures. Vector control is a
key for combating mosquito-borne diseases and the major tool
available for tackling the transmission of Dengue, a disease for
which there is no vaccine nor prophylaxis. There are two main
methods for primary prevention: larval control and adult mosquito
control, depending on the intended target \cite{Natal2002}.
In urban areas \emph{Aedes} mosquitos breed on
water collections in artificial containers such as cans, plastic
cups, used tires, broken bottles, flower pots, etc. Proper solid
waste disposal and improved water storage practices, including
covering containers to prevent access by egg-laying female
mosquitoes, are methods that are encouraged through
community-based programmes \cite{CDC2010}. Active monitoring and
surveillance of the natural mosquito population should accompany
control efforts to determine programme effectiveness.


\section{The mathematical model}
\label{sec:3}

The mathematical model is based on \cite{Dumont2008,Dumont2010}.
The novelty in this paper is the presence of the control parameter
related to adult mosquito insecticide.

The notation used in our mathematical model includes
four epidemiological states for humans:

\begin{quote}
\begin{tabular}{ll}
$S_h(t)$ & susceptible (individuals who can contract the disease)\\
$E_h(t)$ & exposed (individuals who have been infected by the parasite\\
& but are not yet able to transmit to others)\\
$I_h(t)$ & infected (individuals capable of transmitting the disease to others)\\
$R_h(t)$ & resistant (individuals who have acquired immunity)
\end{tabular}
\end{quote}

It is assumed that the total human population $(N_h)$ is constant,
so, $N_h=S_h+E_h+I_h+R_h$.
There are also four other state variables related to the female
mosquitoes (the male mosquitoes are not considered in this study
because they do not bite humans and consequently they do not
influence the dynamics of the disease):

\begin{quote}
\begin{tabular}{ll}
$A_m(t)$& aquatic phase (that includes the egg, larva and pupa stages)\\
$S_m(t)$& susceptible (mosquitoes that are able to contract the disease)\\
$E_m(t)$& exposed (mosquitoes that are infected but are not yet able \\
& to transmit to humans)\\
$I_m(t)$& infected (mosquitoes capable of transmitting the disease to humans)\\
\end{tabular}
\end{quote}

In order to analyze the effects of campaigns to fight the mosquito,
there is also a control variable:

\begin{quote}
\begin{tabular}{ll}
$c(t)$& level of insecticide campaigns\\
\end{tabular}
\end{quote}

Some assumptions in this model:
\begin{itemize}
\item the total human population ($N_h$) is constant; \item there
is no immigration of infected individuals into the human
population; \item the population is homogeneous, which means that
every individual of a compartment is homogenously mixed with the
other individuals; \item the coefficient of transmission of the
disease is fixed and does not vary seasonally; \item both human
and mosquitoes are assumed to be born susceptible; there is no
natural protection; \item for the mosquito there is no resistant
phase, due to its short lifetime.
\end{itemize}

The parameters of the model are:

\begin{quote}
\begin{tabular}{ll}
$N_h$ & total population \\
$B$ & average daily biting (per day)\\
$\beta_{mh}$ & transmission probability from $I_m$ (per bite) \\
$\beta_{hm}$ & transmission probability from $I_h$ (per bite) \\
$1/\mu_{h}$ & average lifespan of humans (in days) \\
$1/\eta_{h}$ & mean viraemic period (in days)\\
$1/\mu_{m}$ & average lifespan of adult mosquitoes (in days) \\
$\mu_{b}$ & number of eggs at each deposit per capita (per day) \\
\end{tabular}
\end{quote}
\begin{quote}
\begin{tabular}{ll}
$\mu_{A}$ & natural mortality of larvae (per day) \\
$\eta_{A}$ & maturation rate from larvae to adult (per day) \\
$1/\eta_{m}$ & extrinsic incubation period (in days)  \\
$1/\nu_{h}$ & intrinsic incubation period (in days) \\
$m$ & female mosquitoes per human \\
$k$ & number of larvae per human \\
$K$ & maximal capacity of larvae
\end{tabular}
\end{quote}

The Dengue epidemic can be modelled by the following nonlinear
time-varying state equations:\\
human Population
\begin{equation}
\label{odehuman}
\begin{tabular}{l}
$\left\{
\begin{array}{l}
\displaystyle\frac{dS_h}{dt}(t)
= \mu_h N_h - \left(B\beta_{mh}\frac{I_m}{N_h}+\mu_h\right) S_h\\
\displaystyle\frac{dE_h}{dt}(t)
= B\beta_{mh}\frac{I_m}{N_h}S_h - (\nu_h + \mu_h )E_h\\
\displaystyle\frac{dI_h}{dt}(t)
= \nu_h E_h -(\eta_h  +\mu_h) I_h\\
\displaystyle\frac{dR_h}{dt}(t)
= \eta_h I_h - \mu_h R_h
\end{array}
\right. $\\
\end{tabular}
\end{equation}
and vector population
\begin{equation}
\label{odevector}
\begin{tabular}{l}
$
\left\{
\begin{array}{l}
\displaystyle\frac{dA_m}{dt}(t)
= \mu_b \left(1-\frac{A_m}{K}\right)(S_m+E_m+I_m)-(\eta_A+\mu_A) A_m\\
\displaystyle\frac{dS_m}{dt}(t)
= -\left(B \beta_{hm}\frac{I_h}{N_h}+\mu_m\right) S_m+\eta_A A_m-c S_m\\
\displaystyle\frac{dE_m}{dt}(t)
= B \beta_{hm}\frac{I_h}{N_h}S_m-(\mu_m + \eta_m) E_m-c E_m\\
\displaystyle\frac{dI_m}{dt}(t)
= \eta_m E_m -\mu_m I_m - c I_m\\
\end{array}
\right. $
\end{tabular}
\end{equation}
with the initial conditions
\begin{equation*}
\label{initial}
\begin{tabular}{llll}
$S_h(0)=S_{h0},$ & $E_h(0)=E_{h0},$ & $I_h(0)=I_{h0},$ &
$R_h(0)=R_{h0},$ \\
$A_m(0)=A_{m0},$ & $S_{m}(0)=S_{m0},$ &
$E_m(0)=E_{m0},$ & $I_m(0)=I_{m0}.$
\end{tabular}
\end{equation*}

Notice that the equation related to the aquatic phase does not have
the control variable $c$, because the adulticide does not produce
effects in this stage of the life of the mosquito. To combat
the larval phase, it would be necessary to use larvicide.
This treatment should be long-lasting and have World Health Organization
clearance for use in drinking water. As we want to study only
a short period of time, this type of treatment has not been considered here.

With the condition $S_h+E_h+I_h+R_h=N_h$, one can, in the example given,
use $R_h=N_h-S_h-E_h-I_h$ and consider an equivalent system for human
population without considering the $R_h$ differential equation.


\section{Basic Reproduction number, equilibrium points and stability}
\label{sec:4}

Let the set
\begin{center}
\begin{tabular}{l}
\small
$\Omega=\{(S_h,E_h,I_h,A_m,S_m,E_m,I_m)\in \mathbb{R}^{7}_{+}:$\\
$ S_h+E_h+I_h\leq N_h, A_m\leq k N_h, S_m+E_m+I_m\leq m N_h \}$
\end{tabular}
\end{center}
be the region of biological interest, that is positively invariant
under the flow induced by the differential system
\eqref{odehuman}-\eqref{odevector} (see Appendix~\ref{appendixAA}).

\begin{thm}
\label{thm:thm1}
Let $\Omega$ be defined as above. Consider also
\begin{center}
\begin{tabular}{l}
$\mathcal{M}=-\left(c (\eta_A + \mu_A)
+ \mu_A \mu_m + \eta_A (-\mu_b + \mu_m)\right)$.
\end{tabular}
\end{center}
The system \eqref{odehuman}-\eqref{odevector}
admits at most two disease free equilibrium points:
\begin{itemize}
\item if $\mathcal{M}\leq 0$, there is a Disease-Free Equilibrium (DFE),
called Trivial Equilibrium, $E_{1}^{*}=\left(N_h,0,0,0,0,0,0\right)$;
\item if $\mathcal{M}> 0$, there is a Biologically Realistic Disease-Free Equilibrium (BRDFE),
$E_{2}^{*}=\left(N_h,0,0,\frac{k N_h \mathcal{M}}{\eta_A\mu_b},
\frac{k N_h \mathcal{M}}{\mu_b \mu_m},0,0\right)$.
\end{itemize}
\end{thm}

\begin{proof}
See Appendix~\ref{appendixA}.
\end{proof}

\begin{remark}
The condition $\mathcal{M}>0$ is
equivalent, by algebraic manipulation, to the condition
$\displaystyle \frac{(\eta_A+\mu_A)(\mu_m+c)}{\mu_b\eta_A}<1$,
where the left-hand side corresponds to the basic offspring number for mosquitos.
Thus, if $M \le 0$, then the mosquito population will collapse and
the only equilibrium for the whole system is the trivial
equilibrium. If $M > 0$, then the mosquito population is
sustainable.
\end{remark}

It is necessary to determine the \emph{basic reproduction number}
of the disease, $\mathcal{R}_0$. This number is very important
from the epidemiologistic point of view. It represents the
expected number of secondary cases produced in a completed
susceptible population, by a typical infected individual during
its entire period of infectiousness \cite{Heffernan2005,
Hethcote2000}. Following \cite{Driessche2002}, we prove:

\begin{thm}
\label{thm:thm2}
If $\mathcal{M}>0$, then the square of the basic
reproduction number associated to
\eqref{odehuman}-\eqref{odevector} is
$\mathcal{R}_{0}^2=\displaystyle\frac{B^2 S_{h} S_{m} \beta_{hm}
\beta_{mh} \eta_m \nu_h }{N_{h}^{2} (\eta_h + \mu_h) (c + \mu_m)
(c + \eta_m + \mu_m) (\mu_h + \nu_h)}$. The equilibrium point
BRDFE is locally asymptotically stable if $\mathcal{R}_{0}<1$ and
unstable if $\mathcal{R}_{0}>1$.
\end{thm}

\begin{proof}
See Appendix~\ref{appendixB}.
\end{proof}

\begin{remark}
In our model, we have two different
populations (human and vector), so the expected basic reproduction
number reflects the infection human-vector and also
vector-human, that is, $\mathcal{R}_{0}^2=R_{hm}\times R_{mh}$ with 
$$
R_{hm} = \frac{B S_{m} \beta_{hm}\nu_h}{N_{h}(\eta_h + \mu_h) (\mu_h + \nu_h)}
$$ 
and 
$$
R_{mh}=\frac{B S_{h} \beta_{mh} \eta_m}{N_{h} (c + \mu_m)(c + \eta_m + \mu_m)}.
$$
The term $B\beta_{hm}\frac{S_{m}}{N_h}$ represents the
product between the transmission probability of the disease from
humans to vectors and the number of susceptible mosquitos per
human; $\frac{1}{\eta_h+\mu_h}$ is related to the human's viraemic
period; and $\frac{\eta_m}{c+\eta_m+\mu_m}$ represents the
proportion of mosquitoes that survive the incubation period.
Analogously, the term $B\beta_{mh}\frac{S_{h}}{N_h}$ is related to
the transmission probability of the disease between mosquitos and
human, in a susceptible population; $\frac{1}{c+\mu_m}$
represents the lifespan of an adult mosquito; and
$\frac{\nu_h}{\mu_h+\nu_h}$ is the proportion of humans who
survive the incubation period.
\end{remark}

When $\mathcal{R}_{0}<1$,
each infected individual produces, on
average, less than one new infected individual, and therefore it
is predictable that the infection will be cleared from the
population. If $\mathcal{R}_{0}>1$, the disease is able to invade
the susceptible population.

\begin{thm}
\label{thm:thm3} If $\mathcal{M}>0$ and $\mathcal{R}_{0}>1$, then
system \eqref{odehuman}-\eqref{odevector} also admits an endemic
equilibrium (EE):
$E_{3}^{*}=\left(S_h^*,E_h^*,I_h^*,A_m^*,S_m^*,E_m^*,I_m^*\right)$,
where
\begin{equation*}
\begin{split}
S_h^* &=
N_h-\displaystyle\frac{(\mu_h+\nu_h)(\mu_h+\eta_h)}{\mu_h\nu_h}I_h^{*},\\
E_h^* &=\displaystyle\frac{\mu_h+\eta_h}{\nu_h}I^{*}_{h},\\
I_h^* &= \frac{\xi}{\chi},\\
\xi &= N_h \mu_h \Bigl[-B^2 k \beta_{hm} \beta_{mh}\nu_h\eta_m
\mathcal{M} + \mu_b \mu_m^2(\eta_m + \mu_m)(\mu_h + \nu_h)(\mu_h + \eta_h)\\
&\quad + c^2\mu_b(\eta_h + \mu_h)(\mu_h + \nu_h)(c + \eta_m + 3 \mu_m)\\
&\quad + c\mu_b\mu_m(\mu_h + \nu_h)(\mu_h(3 \mu_m + 2)
+ \eta_h(2 \eta_m + 3 \mu_m))\Bigr],\\
\chi &= B \beta_{hm} (\eta_h + \mu_h) \Bigl[-\mu_b \mu_h (c +
\mu_m) (c + \eta_m + \mu_m)
- B k \beta_{mh} \eta_m \mathcal{M}\Bigr] (\mu_h + \nu_h),\\
A_m^* &= \displaystyle\frac{\mathcal{M}}{\eta_A \mu_b}k N_h,\\
S_m^* &=\displaystyle\frac{k N_h^2 \mathcal{M}}{\mu_b (c N_h + B I_h^{*} \beta_{hm} + N_h \mu_m)},\\
E_m^* &= \displaystyle\frac{\mu_m+c}{\eta_m}I_m^*,\\
I_m^* &=\displaystyle\frac{B I_h^{*} k N_h \beta_{hm} \eta_m
\mathcal{M}}{\mu_b (c + \mu_m) (c + \eta_m + \mu_m) (c N_h + B
I_h^{*} \beta_{hm} + N_h \mu_m)} .
\end{split}
\end{equation*}
\end{thm}

\begin{proof}
See Appendix~\ref{appendixA}.
\end{proof}

From a biological point of view, it is desirable that humans
and mosquitoes coexist without the disease reaching a level of endemicity.
We claim that proper use of the control $c$ can result in the basic
reproduction number remaining below unity and, therefore, making BRDFE stable.

In order to make effective use of achievable insecticide control,
and simultaneously to explain more easily to the competent authorities its effectiveness,
we assume that $c$ is constant. The goal is to find $c$ such that $\mathcal{R}_{0}^2<1$.


\section{Dengue in Cape Verde}
\label{sec:5}

An unprecedented outbreak was detected in the Cape
Verde archipelago in September 2009 \cite{Who2009reportCV}. 
This is the first report of Dengue virus activity 
in that country. A total of 17 224 cases,
including six deaths, were reported from 18 of the 22
municipalities in Cape Verde by the end of 2009. The municipality
of Praia, on Santiago island, notified the highest number of cases
(13000 cases) followed by S\~{a}o Felipe in Fogo Island (3000 cases).
We used the data for human population related to Cape Verde \cite{CDC2010}.
There have been no cases of Dengue in all the
West African countries near Cape Verde islands since 2000 \cite{Africa21}.
They speak explicitly about mosquitoes coming from Brazil. Also
the information that comes from the Ministry of Health
in the capital of Cape Verde, Praia, confirms
that the insects responsible for Dengue came most probably from Brazil,
transported by means of air transport that perform frequent
connections between Cape Verde and Brazil, as
reported by Jos\'{e} Rosa in the Radio of Cape Verde.
With respect to \emph{Aedes aegypti}, we have thus considered
data from Brazil \cite{Thome2010,Yang2009}.

The simulations were carried out using the following values:
$N_h=480000$, $B=1$, $\beta_{mh}=0.375$, $\beta_{hm}=0.375$,
$\mu_{h}=1/(71\times365)$, $\eta_{h}=1/3$, $\mu_{m}=1/11$,
$\mu_{b}=6$, $\mu_{A}=1/4$, $\eta_A=0.08$, $\eta_m=1/11$,
$\nu_h=1/4$, $m=6$, $k=3$, $K=kN_h$. The initial conditions for
the problem were: $S_{h0}=N_h-E_{h0}-I_{h0}$, $E_{h0}=216$,
$I_{h0}=434$, $R_{h0}=0$, $A_{m0}=kN_h$, $S_{m0}=mN_h$.

Considering nonexistence of control, \textrm{i.e.} $c=0$, the basic
reproduction number for this outbreak in Cape Verde is
approximately $\mathcal{R}_{0}=2.396$, which is in agreement to other
studies of Dengue in other countries \cite{Nishiura2006}.
The control $c$ affects the basic reproduction number,
and our aim is to find a control that puts $\mathcal{R}_{0}$ less than one.

\begin{prop}
Let us consider the parameters listed above and consider $c$ as a constant.
Then $\mathcal{R}_{0}<1$ if and only if $c > 0.156961$.
\end{prop}

Several computational investigations were carried out. 
The software used was \textsf{Scilab} \cite{Campbell2006}.
It is an open source, cross-platform numerical computational
package and a high-level, numerically oriented programming
language. For our problem, we used the routine \texttt{ode} to
solve the set of differential equations. By default, \texttt{ode}
uses the lsoda solver of package ODEPACK. It automatically selects
between the non-stiff predictor-corrector Adams method and the stiff
Backward Differentiation Formula (BDF) method. It uses the non-stiff
method initially and dynamically monitors data in order to decide
which method to use. The graphics were also obtained using this
software, using the command \texttt{plot}.

\begin{figure}[ptbh]
\centering
\begin{minipage}[t]{0.40\linewidth}
\centering
\includegraphics[scale=0.5]{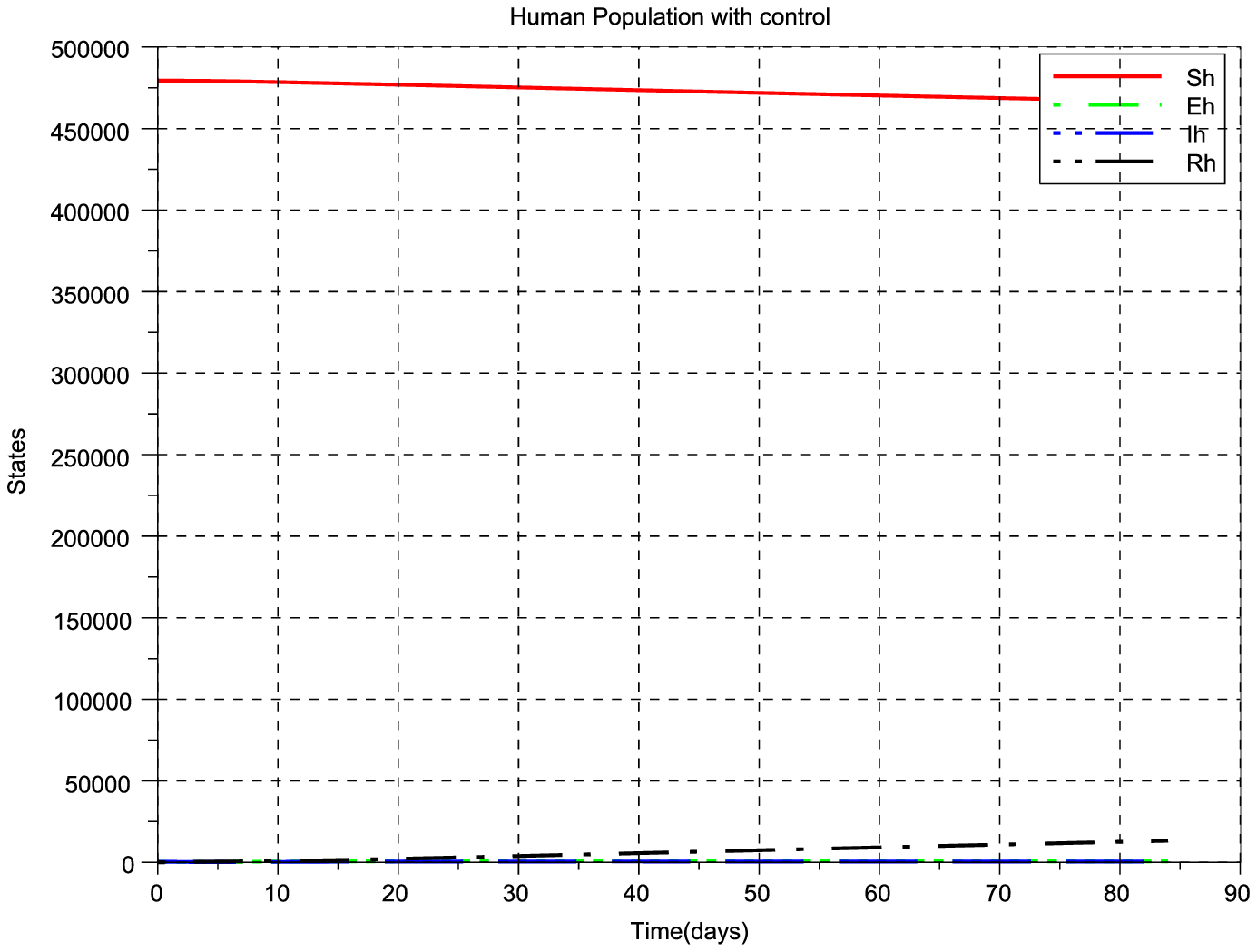}
{\caption{\label{human_control}  Human compartments using control.}}
\end{minipage}\hspace*{\fill}
\begin{minipage}[t]{0.40\linewidth}
\centering
\includegraphics[scale=0.5]{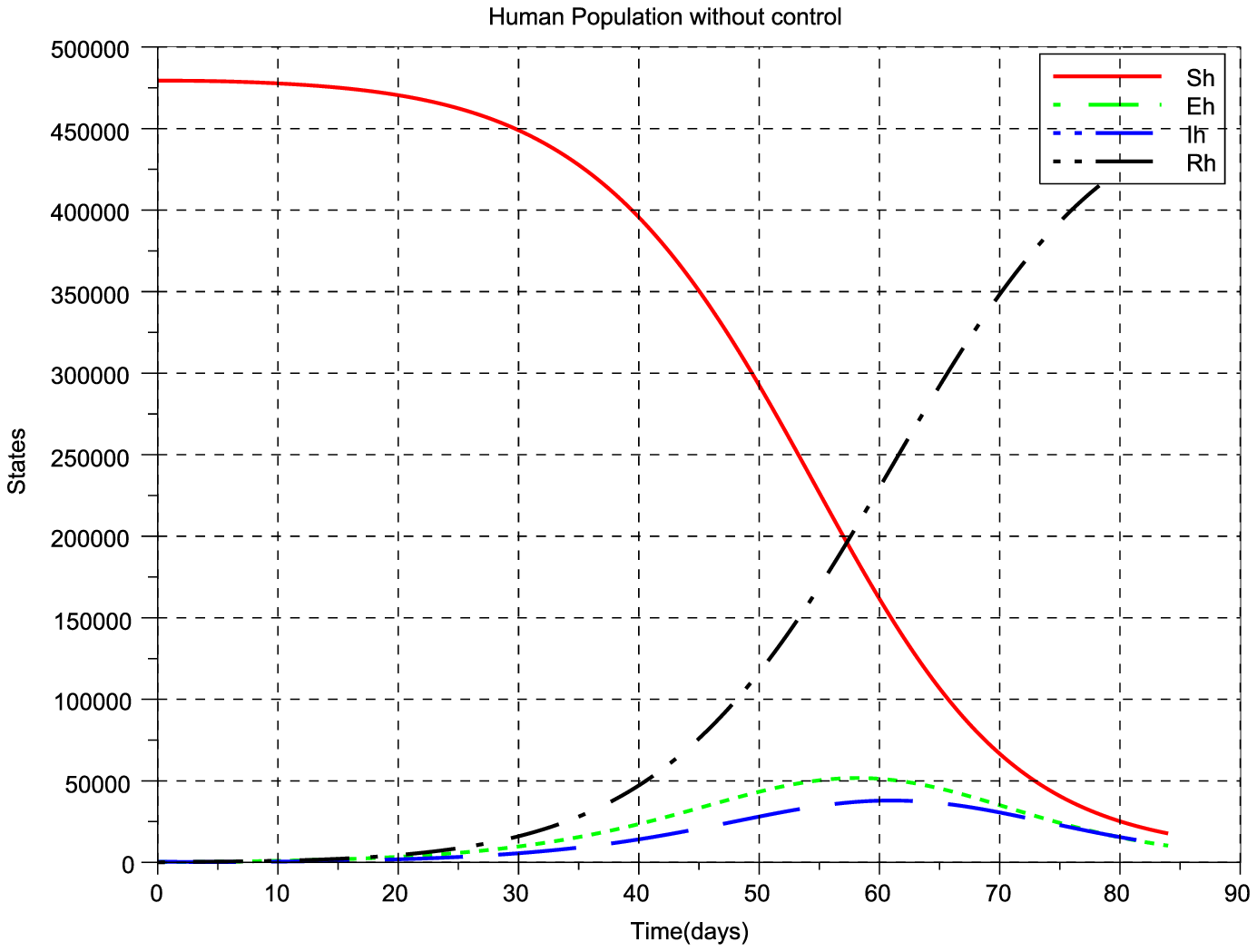}
{\caption{\label{human_nocontrol} \small Human compartments without control.}}
\end{minipage}
\end{figure}

Figures~\ref{human_control} and \ref{human_nocontrol} show the curves
related to human population, with and without control, respectively.
The number of infected persons, even with a small control,
is much less than without any insecticide campaign.

\begin{figure}[ptbh]
\centering
\begin{minipage}[t]{0.40\linewidth}
\centering
\includegraphics[scale=0.50]{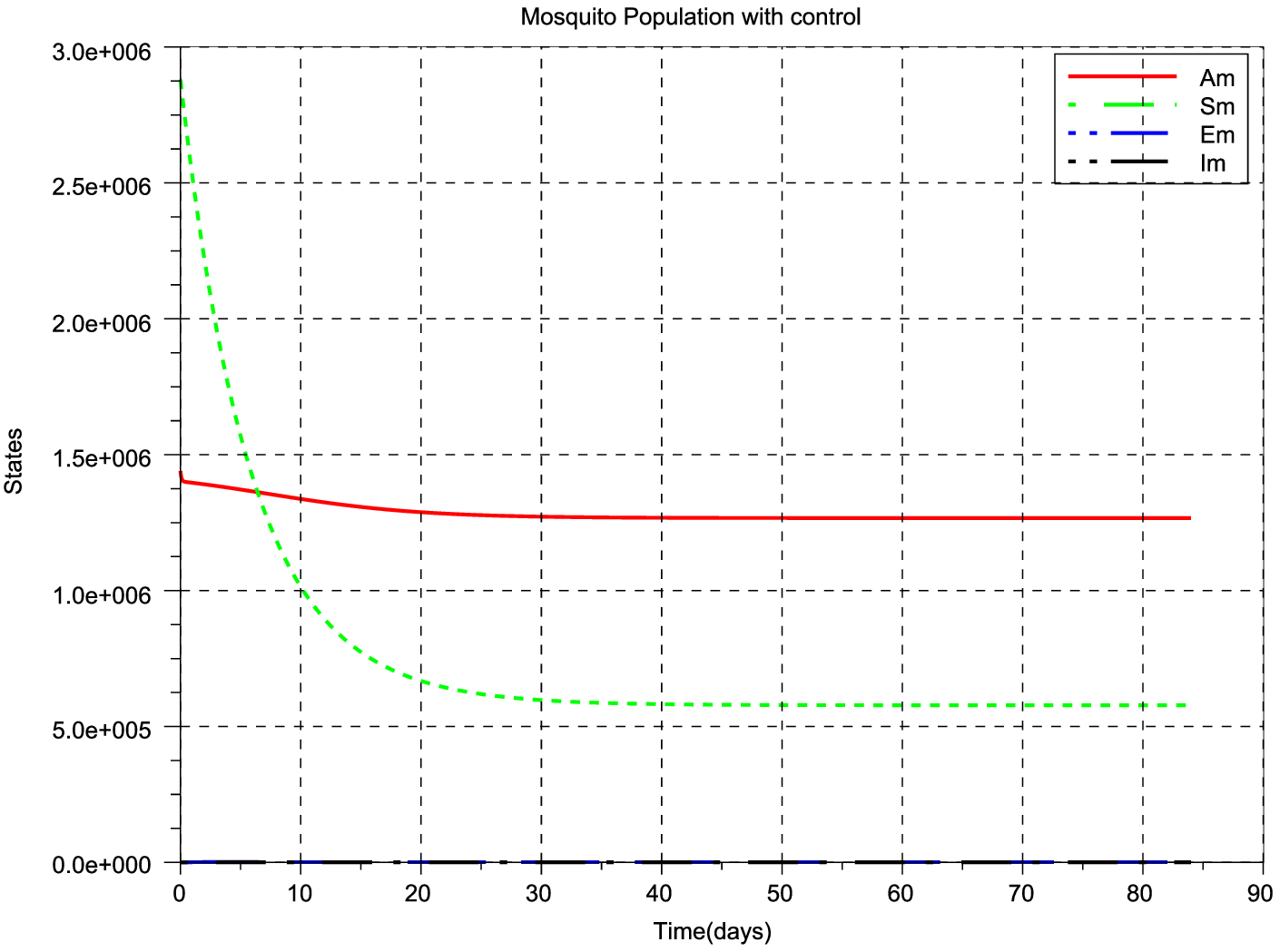}
{\caption{\label{mosquito_control}  Mosquito compartments using control.}}
\end{minipage}\hspace*{\fill}
\begin{minipage}[t]{0.40\linewidth}
\centering
\includegraphics[scale=0.50]{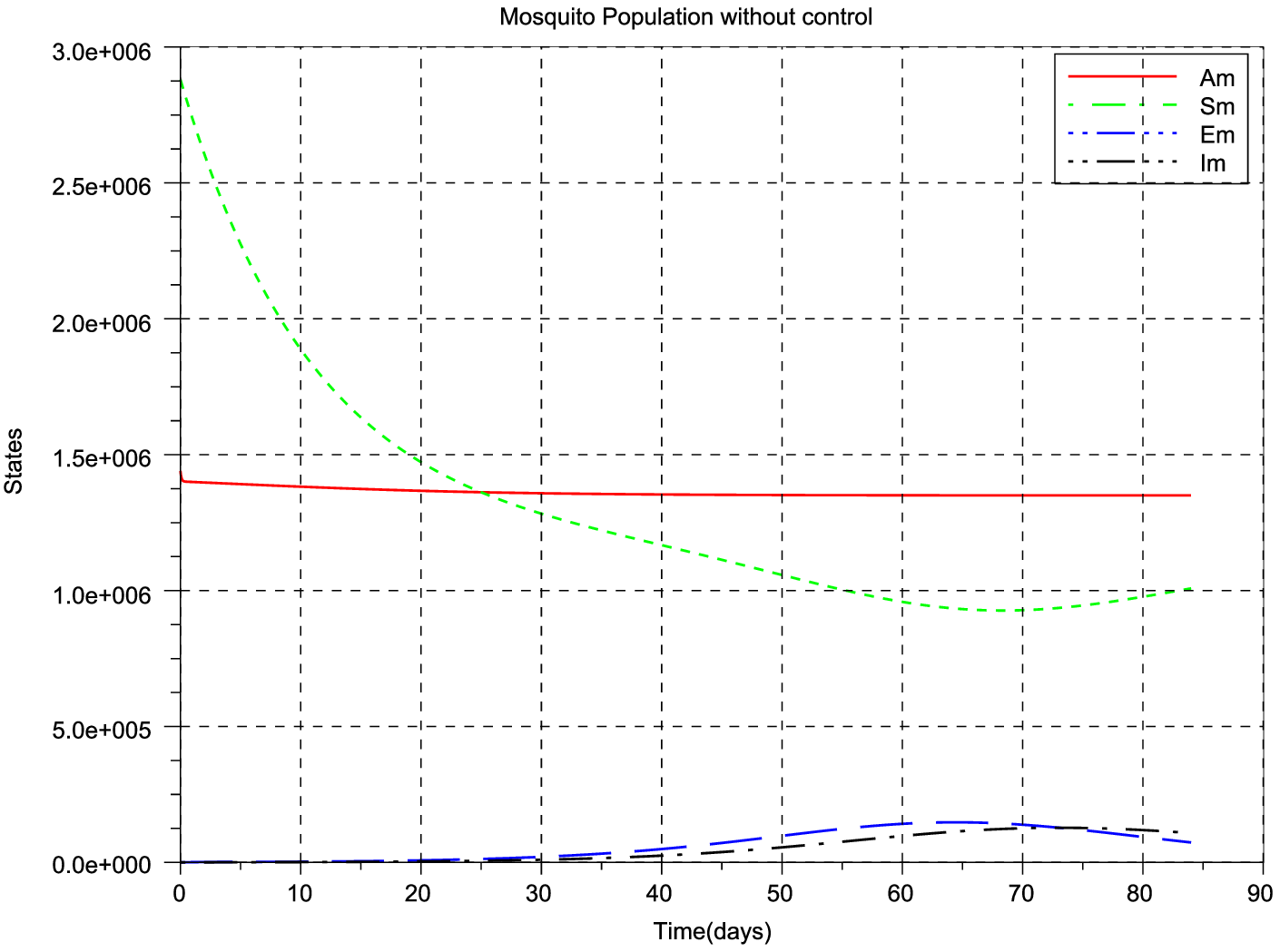}
{\caption{\label{mosquito_nocontrol} \small Mosquito compartments without control.}}
\end{minipage}
\end{figure}

Figures~\ref{mosquito_control} and \ref{mosquito_nocontrol}
show the difference between a region with control and without control.

When the control is applied, the number of infected mosquitoes is close to zero.
Note that the intention is not to eradicate the mosquitoes
but instead the number of infected mosquitoes.


\section{Conclusions and future work}
\label{sec:6}

In this paper, a model based on two populations, humans and mosquitoes,
with insecticide control has been presented.
It has been shown that with a steady insecticide campaign, 
it is possible to reduce the number
of infected humans and mosquitoes and can prevent an outbreak
that could transform an epidemiological episode to an endemic disease.

It has been proved algebraically that, if a constant minimum
level of a control is applied,
it is possible to maintain the basic reproduction number below unity,
guaranteeing the BRDFE. This is corroborated
in another numerical study \cite{Rodrigues2010c}.

In this work, we considered a constant control. In future work,
using a theoretical approach \cite{Rodrigues2010a}, we intend to
find the best function $c(t)$, using optimal control theory.
Instead of finding a constant control, it will be possible
to study other types of control, such as piecewise constant or
even continuous but not constant. Numerical methods can be used to solve the model
\cite{RefereeAsked02,RefereeAsked03,RefereeAsked01}.
Additionally, we could consider another strategy,
a more practical one: due to logistics and health reasons,
it may be more convenient to apply insecticide
periodically and at some specific hours at night.

The rapid increase in mosquito resistance to several
chemical insecticides, and the damage caused by these to the
environment, has resulted in the search for new control
alternatives, such as the use of biological agents. Among the
alternatives available, the use of Bacillus thuringiensis
israelensis (Bti) has been adopted by several countries
\cite{ministeriodasaude2009}. Laboratory testing shows
that Bti has a high larvicide property and its mechanism of action
is based on the production of an endotoxin protein that, when
ingested by the larvae, causes death. So, in future work,
we will also take into account the use of larvicide control.

To ensure the minimization of the outbreaks, educational
programmes that are customized for different levels of health care
and that reflect local capacity should be supported and
implemented widely. The population should also be educated regarding
Dengue in order to minimize the breeding of the
mosquito. Educational campaigns can be included as an extra
control parameter in the model.


\section*{Acknowledgements}

Work partially supported by Portuguese Foundation for Science and Technology (FCT)
through the PhD Grant SFRH/BD/33384/2008 (Rodrigues)
and the R\&D units Algoritmi (Monteiro) and CIDMA (Torres).
The authors are very grateful to three referees,
for valuable remarks and comments, which have
contributed significantly to the quality of the paper.



\appendix

\section{Dynamical properties of Metzler systems}
\label{appendixAA}

System \eqref{odehuman}-\eqref{odevector}
can be rewritten in the following way:
\begin{equation}
\label{odematrix} \displaystyle\frac{dX}{dt}=M(X)X+F,
\end{equation}
where $X=\left(S_h,E_h,I_h,A_m,S_m,E_m,I_m\right)$,
\tiny
\begin{center}
\hspace*{-1.5cm}
\begin{tabular}{c}
$ M(X)=\left( \begin{array}{ccccccc}
- B\beta_{mh}\frac{I_m}{N_h}-\mu_h & 0 & 0 & 0 & 0& 0& 0\\
B\beta_{mh}\frac{I_m}{N_h} & -\nu_h - \mu_h & 0 &0  &0 &0 & 0\\
 0 & \nu_h & -\eta_h  -\mu_h & 0 &0 &0 &0 \\
 0& 0 &0  & -\mu_b\frac{S_m+E_m+I_m}{K}-\mu_m - \eta_m & \mu_b & \mu_b & \mu_b \\
0 & 0 &0  & \eta_A & -B \beta_{hm}\frac{I_h}{N_h}-\mu_m -c &0 &0 \\
 0 & 0 & 0 & 0 &B \beta_{hm}\frac{I_h}{N_h} & -\mu_m - \eta_m -c & 0\\
 0& 0 &0  &0  & 0  &\eta_m & -\mu_m -c\\
\end{array} \right)$,
\end{tabular}
\end{center}
\normalsize
\noindent and $F=\left(\mu_h N_h,0,0,0,0,0,0\right)^{T}$.
As $M(X)$ has all off-diagonal entries nonnegative,
$M(X)$ is a Metzler matrix.

Using the fact that $F\geq 0$, system
\eqref{odematrix} is positively invariant in $\mathbb{R}^{7}_{+}$
\cite{Abate2009}, which means that any trajectory of the system
starting from an initial state in the positive orthant
$\mathbb{R}^{7}_{+}$ remains forever in $\mathbb{R}^{7}_{+}$.


\section{Equilibrium points}
\label{appendixA}

\normalsize

The equilibrium points are reached
when the following equations hold:

\begin{equation}\label{equilibrium}
\begin{tabular}{l}
$\left\{
\begin{array}{l}
\frac{dS_h}{dt}(t) = 0\\
\frac{dE_h}{dt}(t) = 0\\
\frac{dI_h}{dt}(t) = 0\\
\frac{dA_m}{dt}(t) = 0\\
\frac{dS_m}{dt}(t) = 0\\
\frac{dE_m}{dt}(t) = 0\\
\frac{dI_m}{dt}(t) = 0\\
\end{array}
\right. $\\
\end{tabular}
\end{equation}

Using the \textsf{Mathematica} software to solve
system \eqref{equilibrium}, we obtained four solutions.

The first one is known as the \emph{Trivial Equilibrium},
since the mosquitoes do not exist, so there is no disease:
\begin{center}
\begin{tabular}{l}
$E_{1}^{*}=\left(N_h,0,0,0,0,0,0\right)$.
\end{tabular}
\end{center}

In the second one, mosquitoes and humans interact, but there
is only one outbreak of the disease, \textrm{i.e.},
over time the disease goes away without being necessary
to kill all the mosquitoes. We have called this equilibrium point a
\emph{Biologically Realistic Disease-Free Equilibrium} (BRDFE),
since it is a more reasonable situation to find in nature than the previous one:

\begin{center}
\begin{tabular}{lr}
$E_{2}^{*}=$&$\left(N_h,0,0,\displaystyle\frac{k N_h (-\left(c (\eta_A + \mu_A)
+ \mu_A \mu_m + \eta_A (-\mu_b + \mu_m)\right))}{\eta_A\mu_b},\right.$\\
&$\displaystyle\left.\frac{k N_h (-\left(c (\eta_A + \mu_A) + \mu_A \mu_m
+ \eta_A (-\mu_b + \mu_m)\right))}{\mu_b \mu_m},0,0\right)$,
\end{tabular}
\end{center}
which is equivalent to
\begin{center}
\begin{tabular}{l}
$E_{2}^{*}=\left(N_h,0,0,\frac{k N_h \mathcal{M}}{\eta_A\mu_b},
\frac{k N_h \mathcal{M}}{\mu_b \mu_m},0,0\right)$.
\end{tabular}
\end{center}
This is biologically interesting only if $\mathcal{M}$ is greater than 0.

The third solution corresponds to a situation where humans and mosquitoes
live together but the disease persists in the two populations.
So the disease is not anymore an epidemic episode, but transforms
into an endemic one. With some algebraic manipulations, we obtain the following point:

$E_{3}^{*}=\left(S_h^*,E_h^*,I_h^*,A_m^*,S_m^*,E_m^*,I_m^*\right)$,
where

$S_h^*=
N_h-\displaystyle\frac{(\mu_h+\nu_h)(\mu_h+\eta_h)}{\mu_h\nu_h}I_h^{*}$,

$E_h^*=\displaystyle\frac{\mu_h+\eta_h}{\nu_h}I^{*}_{h}$,

$I_h^*=N_h \mu_h (-B^2 k \beta_{hm} \beta_{mh}\nu_h\eta_m
     \mathcal{M} + \mu_b \mu_m^2(\eta_m + \mu_m)(\mu_h + \
\nu_h)(\mu_h + \eta_h) +
    c^2\mu_b(\eta_h + \mu_h)(\mu_h + \nu_h)(c + \eta_m +
       3 \mu_m) +
    c\mu_b\mu_m(\mu_h + \nu_h)(\mu_h(3 \mu_m +
          2) + \eta_h(2 \eta_m +
          3 \mu_m)))/(B \beta_{hm} (\eta_h + \mu_h) (-\mu_b \mu_h \
(c + \mu_m) (c + \eta_m + \mu_m) -
      B k \beta_{mh} \eta_m \mathcal{M}) (\mu_h + \nu_h))$;

$A_m^*=\displaystyle\frac{\mathcal{M}}{\eta_A \mu_b}k N_h$,

$S_m^*=\displaystyle\frac{k N_h^2 \mathcal{M}}{\mu_b (c N_h
+ B I_h^{*} \beta_{hm} + N_h \mu_m)}$,

$E_m^*=\displaystyle\frac{\mu_m+c}{\eta_m}I_m^*$,

$I_m^*=\displaystyle\frac{B I_h^{*} k N_h \beta_{hm}
\eta_m \mathcal{M}}{\mu_b (c + \mu_m) (c + \eta_m + \mu_m)
(c N_h + B I_h^{*} \beta_{hm} + N_h \mu_m)}$.

\medskip

As before, this equilibrium is only biologically interesting if $\mathcal{M}>0$.

With the \textsf{Mathematica} software we obtained a fourth solution.
However, some of its components are negative, which
means that it does not belong to the $\Omega$ set.


\section{The basic reproduction number: proof of Theorem~\ref{thm:thm2}}
\label{appendixB}

The basic reproduction number is calculated in a disease free equilibrium.
In this case we consider the most realistic one, BRDFE.

Following \cite{Driessche2008,Driessche2002}, we consider the vector
$x^{T}=\left(E_h,I_h,E_m,I_m\right)$, which corresponds
to the components related to the progression of the disease.
Thus, the subsystem used is:

\begin{equation}
\label{subsystem}
\begin{tabular}{l}
$
\left\{
\begin{array}{l}
\frac{dE_h}{dt}(t)
= B\beta_{mh}\frac{I_m}{N_h}S_h - (\nu_h + \mu_h )E_h\\
\frac{dI_h}{dt}(t)
= \nu_h E_h -(\eta_h  +\mu_h) I_h\\
\frac{dE_m}{dt}(t)
= B \beta_{hm}\frac{I_h}{N_h}S_m-(\mu_m + \eta_m) E_m-c E_m\\
\frac{dI_m}{dt}(t)
= \eta_m E_m -\mu_m I_m - c I_m .
\end{array}
\right. $\\
\end{tabular}
\end{equation}

This subsystem can be written as partitioned,
$\displaystyle\frac{dx}{dt}=\mathcal{F}(x) - \mathcal{V}(x)$, where $x^{T}
=\left(E_h,I_h,E_m,I_m\right)$, $\mathcal{F}(x)$ represents the components
related to new cases of disease (in this situation in the exposed compartments)
and $\mathcal{V}(x)$ represents the other components.
Thus, the subsystem \eqref{subsystem} can be rewritten as

\medskip

$\mathcal{F}(x)=\left(
                  \begin{array}{c}
                    B\beta_{mh}\frac{I_m}{N_h}S_h \\
                    0 \\
                    B\beta_{hm}\frac{I_h}{N_h}S_m \\
                    0 \\
                  \end{array}
                \right)
$
and
$\mathcal{V}(x)=\left(
                  \begin{array}{c}
                    (\nu_h+\mu_h)E_h\\
                    -\nu_h E_h + (\eta_h+\mu_h)I_h \\
                    (\mu_m+\eta_m+c)E_m\\
                    -\eta_m E_m +(\mu_m+c)I_m \\
                  \end{array}
                \right).
$

\medskip

Let us considerer the Jacobian matrices associated
with $\mathcal{F}$ and $\mathcal{V}$:

\begin{center}
\begin{tabular}{l}
$J_{\mathcal{F}(x)}=\left(
                  \begin{array}{cccc}
                    0 & 0 & 0 & B\beta_{mh}\frac{S_h}{N_h} \\
                    0 & 0 & 0 & 0\\
                    0 & B\beta_{hm}\frac{S_m}{N_h} & 0 & 0\\
                    0 & 0 & 0 & 0\\
                  \end{array}
                \right)
$
\end{tabular}
\end{center}

and
\begin{center}
\begin{tabular}{l}
$J_{\mathcal{V}(x)}=\left(
                  \begin{array}{cccc}
                    \nu_h+\mu_h & 0 & 0 & 0 \\
                    -\nu_h & \eta_h+\mu_h & 0 & 0 \\
                    0 & 0 & \mu_m+\eta_m+c & 0\\
                    0 & 0 & -\eta_m & \mu_m+c\\
                  \end{array}
                \right).
$
\end{tabular}
\end{center}

According to \cite{Driessche2002}, the basic reproduction number is
$\mathcal{R}_{0}=\rho(J_{\mathcal{F}(x_{0})}J_{\mathcal{V}^{-1}(x_{0})})$,
where $x_{0}$ is a disease free equilibrium (BRDFE)
and $\rho(A)$ defines the spectral radius
of a matrix $A$. Using \textsf{Mathematica}

\begin{center}
\begin{tabular}{l}
$\mathcal{R}_{0}^2=\displaystyle\frac{B^2 k \beta_{hm} \beta_{mh}
\eta_m \nu_h\mathcal{M} }{\mu_b (\eta_h + \
\mu_h) \mu_m (c + \mu_m) (c + \eta_m + \mu_m) (\mu_h + \nu_h)}$,
\end{tabular}
\end{center}
and we obtain the value for the threshold parameter, with $\mathcal{M}>0$.


\end{document}